\documentclass[10pt,twoside]{article}
\usepackage{graphicx}
\usepackage{amsmath}
\usepackage{Latex-document}

\title{\bf The Third Approach to the History of \vskip -2mm Mathematics in China\vskip 6mm}

\author{Anjing Qu\vspace*{-0.5cm}\thanks{Northwest University,
Department of Mathematics, Xi'an 710069, China.
 E-mail: qaj@sein.sxgb.com.cn}}

\date{\vspace{-8mm}}

\begin{document}

\maketitle

\thispagestyle{first} \setcounter{page}{947}

\begin{abstract}

\vskip 3mm

The first approach to the history of mathematics in China led by Li Yan (1892--1963) and Qian Baocong (1892--1974)
featured discovering {\it what} mathematics had been done in China's past. From the 1970s on, Wu Wen-tsun and
others shifted this research paradigm to one of recovering {\it how} mathematics was done in ancient China. Both
approaches, however, focus on the same problem, that is mathematics in history. The theme of the third approach is
supposed to be {\it why} mathematics was done. Combining this approach with the former two, the research paradigm
will be improved from one of mathematics in history to that of the history of mathematics.

\vskip 4.5mm

\noindent {\bf 2000 Mathematics Subject Classification:}  01A25.

\noindent {\bf Keywords and Phrases:} Chinese mathematics, Research paradigm, Interpolation, Numerical method,
Scientific tradition.
\end{abstract}

\vskip 12mm

\section{Introduction} \label{section 1}\setzero

\vskip-5mm \hspace{5mm}

Since the beginning of the last century hundreds of scholars have devoted themselves to the discipline of the
history of mathematics in China. Their research has not only thrown light on the various features of traditional
Chinese mathematics, but also has led to a better understanding of the diversity of the mathematical sciences.

This research has led to problems, however.  Some mathematicians complained that most Chinese historians of
mathematics limited their research to ancient China, while it has seemed to other scholars that fresh avenues into
the history of traditional Chinese mathematics may, to some extent, have been exhausted. Fewer and fewer young
scholars are attracted to the field, and even for some senior historians of mathematics it has been difficult to
find exciting new topics to work on.

A similar perception prevailed once before in the 1970s, when many Chinese historians of mathematics had become
discouraged about the future. It came as something of a surprise that, soon after some scholars left to work on
other subjects, a fresh upsurge of interest led by Wu Wen-tsun appeared, one that lasted until a few years ago.

The aim of this article is twofold. First, research paradigms adopted by Chinese historians of mathematics will be
outlined, followed by a discussion of a new way to approach the history of mathematics in China.

\section{Discovery: the first approach} \label{section 2} \setzero

\vskip-5mm \hspace{5mm}

The question of {\it what} mathematical science, if any, existed in ancient China was first raised by scholars at
the beginning of the 20th century, and was the motivation for those who followed turned their attention to the
history of pre-modern Chinese mathematics. There is no doubt that Li Yan (1892--1963) and Qian Baocong
(1892--1974), the founders of the subject of the history of mathematics in China, deserve to be named as the
representative figures of this movement for discovering {\it what} mathematics was done in ancient China. Two
examples, more or less related to them, can be taken in order to show what the word {\it discovery} has meant in
studies in the history of mathematics in China.

\subsection{ Interpolation}

\vskip-5mm \hspace{5mm}

In Yi-xing's {\it Dayan li} (a calendar-making system of 724 AD), a function $f(x)$ designed for calculating the
solar equation of center is found as follows:
\begin{eqnarray}
  f(x)=\frac{x}{n_1}\times \triangle_1+(1-\frac{x}{n_1})\times \frac{x}{2n_2}\times
  \triangle^{2} \label{1}
\end{eqnarray}
where $0\le x<n_{1}$. While a tropical year is broken into 24
parts $(qi)$, $n_{1}$ and $n_{2}$ are the lengths in days of two
consecutive $qi$. $\triangle_{1}$ and $\triangle_{2}$ are the
deviations in $du$ (1 $du=360^{\circ}/365.25$) from the mean solar
motion to its true one on the intervals $n_{1}$ and $n_{2}$
respectively.
$\triangle^{2}=\dfrac{2n_{1}n_{2}}{n_{1}+n_{2}}(\dfrac{\triangle_{1}}{n_{1}}-
\dfrac{\triangle_{2}}{n_{2}})$.Suppose $n=n_{1}=n_{2}$, the
relation yields to $\triangle^{2}=\triangle_{1}-\triangle_{2}$.
This special case, formula (1), is found in Liu Zhuo's {\it
Huangji li} (a calendar-making system of 600 AD).

It was Yabuuti Kiyosi (1906--2000) who pointed out for the first time that Liu Zhuo's formula is a quadratic
interpolation of equal interval, while Yi-xing's formula is that of unequal interval. Both of them are equivalent
to Gauss's interpolation [1].

A decade after Yabuuti's discovery, a more detailed investigation of this topic was made by Li Yan, who
demonstrated that Liu Zhuo's quadratic interpolation occupied a leading position among various numerical methods
in ancient Chinese mathematical astronomy [2].

Assume that the values of a real function $f(x)$ are given at each of $n + 1$ distinct real values
$x_{k}:f_{k}(k=0,1,2,\cdots,n)$. The method of finding the values $f(x)$ at $x$ by using these values
$f_{k}=f(x_{k})$ is called interpolation. Formulae of interpolation could be stated in many ways, for instance the
formulae of Lagrange, Aitken, Newton, Gauss, Stirling, Bessel and Everett, depending on the method you make use of
to construct it. In cases where the points of interpolation one chooses are the same, their interpolation
functions can be transformed to each other no matter by which kind of interpolation the functions are constructed.
As for formula (1), it is easy to verify that
$$f(0)=0,\quad f(n_{1})=\triangle_{1},\quad f(n_{1}+n_{2})=\triangle_{1}+\triangle_{2}.$$
Clearly, $x=0,n_{1}$, and $n_{1}+n_{2}$ are three points of interpolation of the function $f(x)$.This result shows
that formula (1) is a quadratic interpolation function.

Yabuuti said that formula (1) is equivalent to Gauss's interpolation, while Li Yan considered it to be equivalent
to Newton's. Actually, it is neither Gauss's nor Newton's.

What Yabuuti and Li Yan did, as such, was to reveal the fact that formula (1) is a quadratic interpolation
function. However, the question of which kind of interpolation Chinese mathematicians made use of to construct
formula (1) was left open.

\subsection{Remote measurement}

\vskip-5mm \hspace{5mm}

In his {\it Haidao suanjing} (Sea Island Computational Canon, 263 AD), Liu Hui designed nine questions in order to
demonstrate the problems of remote measurement. The first of them is a problem concerning how to measure the
altitude of an island with two gnomons.

Let $HI$ in Fig.1 be the altitude of an island. $AB$ and $CD$ are two gnomons of equal height. $AE$ and $CF$ are
the length of shadows of $AB$ and $CD$, respectively. Suppose the distance $(AC)$ between the two gnomons is
known, a formula for measuring the island's altitude is found in the {\it Haidao suanjing} as follows:
\begin{eqnarray}
  HI=AB+\frac{AB \times AC}{CF-AE} \label{2}
\end{eqnarray}

Formula (2) is essential in remote measurement. Other questions in Liu's book are much more complex than this. As
many as four gnomons are used in some of them. It is said that there were diagrams drawn by Liu Hui for these
questions, but they no longer exist today.

The question of the island's altitude is also known as that of solar altitude.  Formula (2) was used to measure
the altitude of the sun in the {\it Zhoubi suanjing} (Zhou Dynasty Canon of Gnomonic Computations, first century
BC). A diagram for proving formula (2) is also found in this work. Unfortunately, due to transmission of the text
over time, it has been distorted beyond recognition.

For deducing formula (2) or showing its correctness, Qian Baocong added a line $DG$ parallel to the line $HE$ in
Fig.1.[3] What Qian did was not unusual for historians of mathematics at that time.

\centerline{\includegraphics[width=12cm]{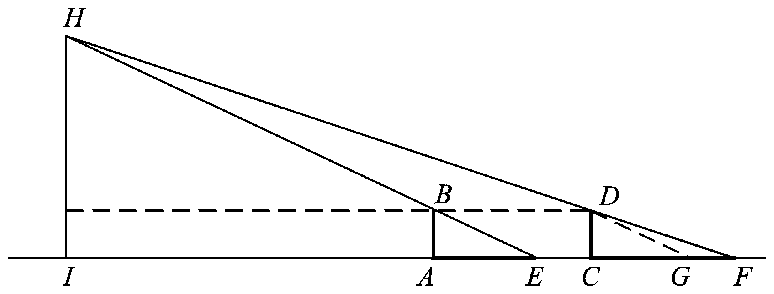}} \vskip 2mm
\centerline{Fig. 1\quad Altitude of island } \vskip 2mm

However, a new movement in this area of research was heralded in
the 1970s when a mathematician who adopted a methodology
different to that of Qian's appeared.

\section{Recovery: the second approach} \label{section 3} \setzero

\vskip-5mm \hspace{5mm}

It was Wu Wen-tsun who led studies in the history of mathematics in China on to a second phase, that of recovering
{\it how} mathematics was done in ancient China. Wu explained his approach by criticizing Qian's study on formula
(2). Since the parallel line $DG$ drawn by Qian in Fig.1 is groundless in traditional Chinese mathematics, his
proof of formula (2), as Wu pointed out, should be regarded as a ``wrong proof'' from the viewpoint of the history
of mathematics.

Wu emphasized that demonstrating the correctness of ancient mathematics with modern mathematical notions is by no
means the sole purpose of the history of mathematics, and that historians of mathematics should pay more attention
to recovering {\it how} mathematics was actually done in history. He said:

{\it ``Two basic principles of such studies will be strictly observed, viz.:

``P1. All conclusions drawn should be based on original texts fortunately preserved up to the present time.

``P2. All conclusions drawn should be based on reasoning in the manner of our ancestors in making use of knowledge
and in utilizing auxiliary tools and methods available only at that ancient time.''} [4]

Wu, therefore, named his approach recovering mathematical procedure with its original thought.

\begin{table}[th]
\centerline{\includegraphics[width=12cm]{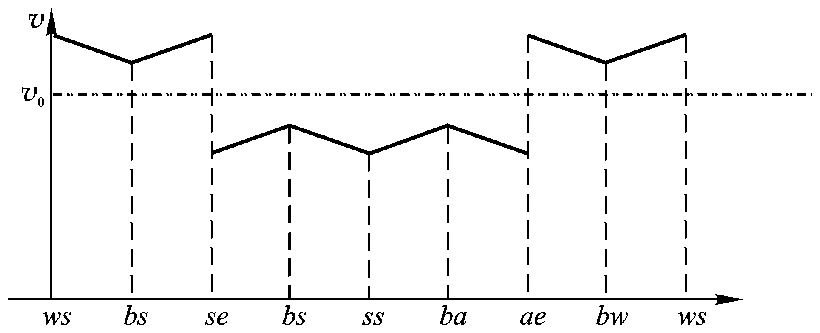}}\vskip 2mm
\centerline{Fig. 2\quad Pattern of solar speed in the {\it
Huangji li} (600 AD)} \vskip 2mm

\centerline{\includegraphics[width=12cm]{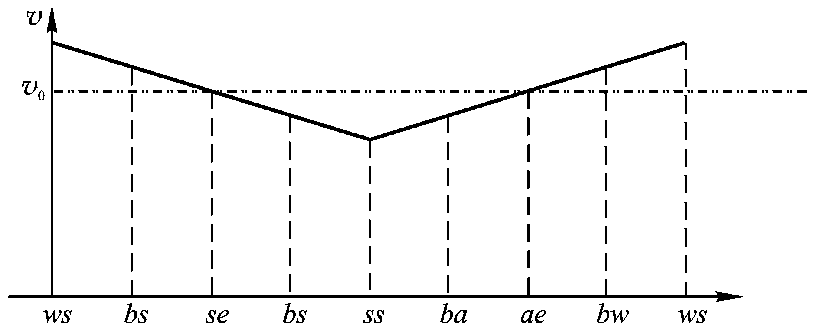}}\vskip 2mm \centerline{Fig. 3\quad Pattern of solar speed in
the {\it Dayan li} (724 AD)}
\end{table}

Let us also take interpolation as an example to show how to recover the history of mathematics in China. How
formula (1) was constructed was exactly the topic that Wu's approach focused on. In order to answer this question,
the historical background from which the problem of interpolation arose needed to be addressed.

It is well known that formula (1) was invented to solve the problem of irregular solar motion. In the {\it Chapter
of Calendar-making} of the {\it Tang History}, the evolution of solar theory up to that time was described by the
Buddhist monk Yi-xing (683--727) as follows:

In about 560 AD, Zhang Zixin found that solar motion was irregular. From Yi-xing's account, we understand that
solar motion in the minds of Liu Zhuo (600 AD) and Yi-xing (724 AD) can be expressed as the patterns in Fig.2 and
Fig.3, respectively, where the dotted line $v_0$ represents the mean velocity of the sun.

Liu's interpolation function is different from Yi-xing's in the division of a year. Liu divided a tropical year
into 24 parts of equal length, while Yi-xing divided the ecliptic into 24 parts of 15$^\circ$ each. The small
interval was named mean {\it qi} ({\it ping qi}) in Liu's division, and ture {\it qi} ({\it ding qi}) in
Yi-xing's. Since solar motion is irregular, Yi-xing's division is unequally spaced in terms of time. In order to
deal with the deviation from the mean motion to the true motion of the sun, interpolation function (1) was
constructed by Liu Zhuo and Yi-xing at each {\it qi}.

The period between the winter solstice ($ws$ in Fig.2--4) and the
beginning of spring ($bs$ in Fig.2--4) consists of 3 {\it qi}. Let
us take this period as an example to demonstrate how formula (1)
was constructed by Liu Zhuo and Yi-xing.

In Fig.4, suppose $OM = n_1$, $MN = n_2$, area of $OBCM = \Delta_1$, area of $MEFN = \Delta_2$. $\Delta_1$ and
$\Delta_2$ given by observation are deviations of the true motion of the sun from its mean motion on the interval
$OM$ and $MN$, respectively. The doted line $OMN$ stands for the velocity of the mean sun, while the step-like
lines $BC$, $EF$, and so on stand for the mean velocities of the true sun on each {\it qi}.

The idea that Liu Zhou and Yi-xing conceived was how to change the pattern of solar motion from a step pattern to
a continuous straight line. The observed data $\Delta_1$ and $\Delta_2$ of the two consecutive {\it qi} were used
to construct the interpolation function on the former {\it qi} so that the slanted lines as a whole could be
continuous lines as far as possible. This is a change from a linear to a parabolic interpolation.

\begin{center}

\includegraphics[width=6cm]{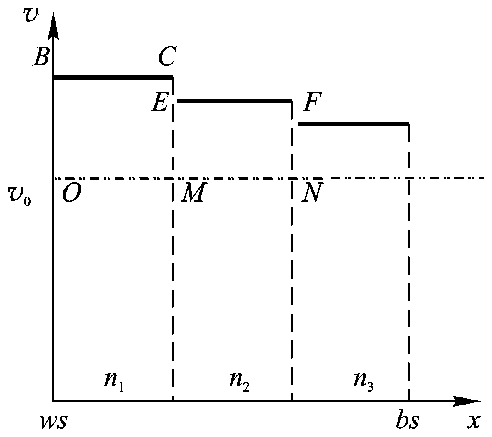} \quad \includegraphics[width=6cm]{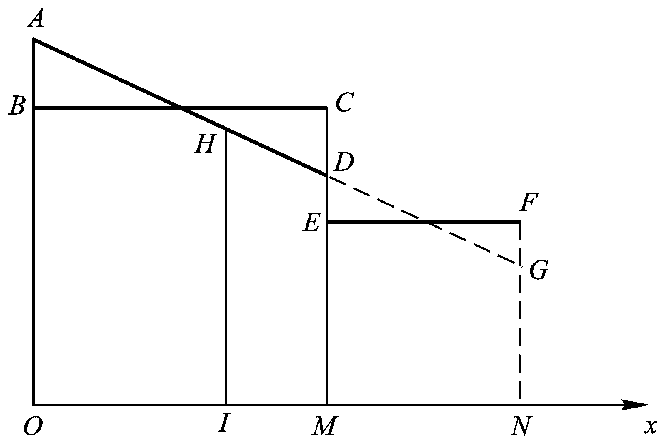}

Fig. 4 \quad Linear interpolation \hspace{1cm} Fig. 5 \quad Quadratic interpolation
\end{center}

As in Fig.5, draw a slanted line $AG$ that crosses the midpoints of lines $BC$ and $EF$, and intersects line $CM$
at $D$. Then $AB = DC$, $DE = FG$. Let $x = OI$, then
$$ f(x) = \mbox{area of trapezoid $AOIH$} $$
is the Chinese interpolation formula (1), $0 \le x < n_1$.

In order to construct the function $f(x)$, first of all Liu and Yi-xing let the solar speed be changed as an
arithmetic sequence from day to day at the interval $OM$. Then, the arithmetic progression is summed up. The
result of the summation is a parabolic function with a variable of days after the initial time at the interval
$OM$. It is the quadratic interpolation function $f(x)$. From the viewpoint of the technique of constructing
algorithms, Yi-xing's unequal spaced interpolation is nothing more than Liu Zhuo's equal intervals interpolation.

Since the tropical year is divided into 24 small intervals ({\it qi}), and the parabolic interpolation function
was constructed on each interval based on the observed data of 24 {\it qi}, we call a whole set of the 24
functions of a tropical year {\it piecewise parabolic interpolation}.

\section{Mathematics in history: original research} \label{section 4} \setzero

\vskip-5mm \hspace{5mm}

During the period of the first movement in the history of Chinese mathematics, {\it discovery} meant to find out
{\it what} mathematics was done in history. Scholars made their discoveries directly from original historical
materials. In the next movement, research discoveries were extended to {\it recoveries}. The attention of scholars
was directed to the question of {\it how} mathematics was done in history.  Recovery work in the history of
mathematics is that of rational reconstruction, usually based on indirect historical materials. Recovery,
therefore, can be regarded as a kind of indirect discovery.

The paradigm of studies in the history of mathematics in China took shape in the following way: only {\it
discovery}, either direct or indirect, was regarded as original research. Once such a paradigm became the norm in
the Chinese history of mathematics community, the model became set for such research.

Obviously, most results gained from the first movement were transformed into the problems to be solved in the
second. Problems raised in the first movement became conjectures, while the use of secondary historical materials
to ``prove'' these conjectures was the main trend in the second movement.

This was the research paradigm that historians of mathematics in China followed during the past century. Studies
in the history of mathematics were valued according to the viewpoint of this paradigm. It is somewhat equivalent
to that of pure mathematics. What they did was to discover or recover mathematics in history.

This picture can help us to clarify the following questions that may puzzle those who are not involved in this
field.

First of all, Wu's movement was the consequence of the fact that the research paradigm shifted from {\it discovery
to recovery} after the 1970s. This change was so important that it offered plenty of topics for research to
historians of mathematics in the last quarter of the last century.

Secondly, as we have mentioned before, original research in the history of mathematics meant research in original
historical materials. For most Chinese scholars, unfortunately, the only original historical materials that they
had access to were texts in Chinese. Historical works in Western mathematics or modern mathematics might interest
mathematicians or lay readers, and might be widely welcomed, but neither discovery nor recovery could be expected.
That meant that these attempts could not be considered as original research. This is the reason why most Chinese
history of mathematics was limited to research on ancient China.

\centerline{\includegraphics[width=12cm]{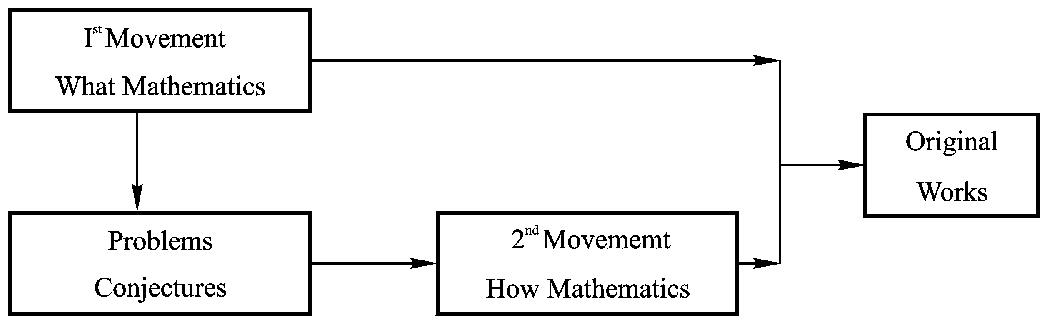}}\vskip 2mm
\centerline{Fig. 6\quad Studies in the history of mathematics in
China: past and present} \vskip 2mm

\section{Why mathematics: the third approach} \label{section 5} \setzero

\vskip-5mm \hspace{5mm}

Discovering {\it what} mathematics was done is always the basic approach to the history of mathematics. During the
period of Li and Qian's movement, it had been taken as the only way to do the history of mathematics in China.
Recovering {\it how} mathematics was done was never seriously accepted as original research until Wu Wen-tsun
broke through Li and Qian's paradigm in the 1970s. It helped Chinese historians of mathematics to negotiate
through a crisis, and fruitful results have come from this extension of the concept of {\it original research} in
the history of mathematics in China. The model depicted in Fig.6 clarifies the work that Chinese historians of
mathematics have done.

The problem now is why are we faced with a crisis once again? How can we step out of such a difficult situation
and develop an optimistic future for our field in China? These are vital issues which I have taken upon myself to
tackle in the last part of this article.

The reason why the field faces a crisis once again is because the concept of {\it original work}, either {\it
discovery} or {\it recovery}, has restricted research on {\it mathematics in history}. Since the original
historical materials that our Chinese scholars can access are limited, they are certain to be exhausted sooner or
later.

The answer to achieving our goal of renewed vigor in the field is that the concept of {\it original work} in the
history of mathematics in China has to be extended further. The old research paradigm, therefore, must be improved
once again.

Mathematical ideas are the main object of study in the history of mathematics. Hence the history of mathematics,
to a great extent, is the history of mathematical thought. When we look at mathematics in history from a
historical perspective, an important aspect is often neglected, that is, {\it why} did mathematics play a part in
history? As I have already emphasized, in the past century, two movements dominated studies in the history of
mathematics in China namely:

The first approach: {\it what mathematics} was done.

The second approach: {\it how mathematics} was done.

It is easy to follow these approaches up with:

A third approach: {\it why mathematics} was done.

Once the three approaches have merged into a single whole, the research paradigm will shift from mathematics in
history to the history of mathematics. {\it Original research} in the history of mathematics will also be extended
in scope.

Actually, {\it why mathematics} has been regarded as the main purpose of the history of mathematics by leading
mathematicians for some time. Andr\'e Weil at ICM 1978, for instance, presented a plenary speech about ``{\it for
whom does one write a general history}'' of mathematics. At the end of his presentation, he said  ``{\it thus my
original question `why mathematical history?' finally reduces itself to the question `why mathematics?' which
fortunately I do not feel called upon to answer.}'' [5] Wu Wen-tsun himself also sometimes moved beyond how
mathematics in ancient China [6].

\section{Why a practical tradition in ancient China} \label{section 6} \setzero

\vskip-5mm \hspace{5mm}

It is often said that, compared with Greek mathematics, Chinese mathematics was characterized by a practical
tradition. Many scholars hold that this tradition is the fatal weakness of Chinese mathematical science, one that
prevented it from developing into modern science. Some historians of mathematics have argued that certain
fundamental factors of the Greek theoretical tradition, such as {\it proof} and {\it principle}, can also be
detected in the {\it Nine Chapters of Arithmetic} (1 century BC) and Liu Hui's {\it Annotations} (263 AD).
However, it seems to us that many people are still not convinced.

A scientific tradition represents the principal aspects of scientific methodology, spirit, and style. When we
speak of the theoretical tradition of Greek science, we do not mean that there was no applied science at all. What
we mean is that compared with this theoretical tradition, the practical tradition was of trivial importance in the
development of Greek science.

For a better understanding of the value of Chinese mathematics from a historical perspective, we need to know
about the issue of {\it why mathematics} was done in ancient civilizations, and thus why there was a practical
tradition in ancient China is a contact point for this.

In the long history of the Chinese empire, mathematical astronomy was the only subject of the exact sciences that
attracted great attention from rulers. In every dynasty, the royal observatory was an indispensable part of the
state. Three kinds of expert --- mathematicians, astronomers and astrologers --- were employed as professional
scientists by the emperor. Those who were called mathematicians took charge of establishing the algorithms of the
calendar-making systems. Most mathematicians were trained as calendar-makers. Mathematics was thus highly
developed for mathematical astronomy besides more general applications, in such areas as indeterminate problems,
numerical solutions of algebraic equations, polynomial interpolation and series summation.

Calendar-makers were required to maintain a high degree of precision in prediction. Ceaseless efforts to improve
numerical methods were made in order to guarantee that the algorithm could satisfy the precision required for
astronomical observation.[7] It was neither necessary nor possible that a geometric model could replace numerical
method, which occupied the principal position in Chinese calendar-making systems. The reason for this was that
only the numerical method could satisfy the ruler's requirements, that is high accuracy in prediction and
computation. As a result numerical analyses won favor over cosmic or geometric model building. As a subject
closely related to numerical method, algebra, rather than geometry, became the most developed field of mathematics
in ancient China.[8]

Science in ancient China was intended primarily to solve concrete problems, such as determining planetary
positions. The function of {\it explaining} natural phenomena never dominated its scientific tradition. What
Chinese scientists really cared about was how to solve the problems they faced as accurately as possible.

This is the reason why the practical tradition was chosen in ancient China. From the above description,
definitions for the practical and theoretical traditions maybe drawn as follows:

In the practical tradition, science serves to solve concrete problems. Theory is judged by its accuracy of
computation. Scientific progress follows the advancement of observation. A theoretical model is always improved to
meet the precision requirement step by step.

In the theoretical tradition, on the other hand, science serves to explain natural phenomena. Theory is judged by
its function in the explanation. Observation is employed to verify the correctness of the theoretical hypothesis.
The old model is always replaced if the new one is more reasonable for the explanation of natural phenomena.

These two traditions differ mainly in their starting and ending points:

In the practical tradition, a model is built up from observation to solve concrete problems. For a more accurate
prediction, uninterrupted efforts are made to explore unknown factors, and the related numerical analyses are
improved. The more accurate the theory is, the closer the model is to the truth.

In the theoretical tradition, on the other hand, a model is built up from hypotheses that account for natural
phenomena. For a better understanding of natural phenomena, the model is revised from time to time on the basis of
a new hypothesis. The closer the model is to the truth, the more accurate is the theory.

The attitude to algebraic equations of these two traditions provides a typical example to show us their
differences. In the theoretical tradition, mathematicians paid attention to the root formula of the equation. The
numerical solution of the equation was seldom taken much notice of by them. The reason is that no matter how
effective the numerical method is, its solution is usually an approximate result that does not help to {\it
explain the phenomenon} of the equation.

On the contrary, the root formula is never more important than a numerical solution in the practical tradition.
The reason is that even if one could have the exact root from the formula, one has to extract from it a concrete
value for application. Hence it seems to mathematicians in the practical tradition that the numerical solution is
sufficient. In ancient China, the root formula of equations was not an important subject, although they did know
the formula of second-degree equations.

It is believed that modern science, to a large extent, benefits from the heritage of Greek science. Nevertheless,
it is hard to say that the theoretical tradition dominates the development of modern science in all aspects. In
fact, the practical tradition also plays an important role.

It is obvious that numerical analyses are more frequently used than theoretical hypotheses in modern science. The
task for modern scientists is not only to account for natural phenomena, but also to solve concrete problems. The
research results arising from the search for solutions to scientific problems have led in two directions: those
that are concerned with finding general theorems concerning the problems, and those that are searching for good
approximations for solutions. Both explaining natural phenomena and solving concrete problems are the goals that
modern scientists strive for. Observations have occupied a substantial position in the development of modern
science.

Generally speaking, science in ancient civilizations was often characterized by a distinctive tradition, either
the theoretical tradition as in Greece, or the practical tradition as in China. These traditions tended to develop
along their own lines. However, the situation in modern times has never been so simple. The diversity of modern
science features the blending of the two traditions. It develops in a dualistic mode.

\section{Conclusion} \label{section 7} \setzero

\vskip-5mm \hspace{5mm}

The paradigm of the history of mathematics in China has directed the attention of researchers to focus on {\it
mathematics in history}, in particular in ancient China. It is certain that the {\it what} and {\it how
mathematics} are two approaches for historians of mathematics that will always remain valid. If history continues,
they continue.

However, just as Li and Qian's what mathematics approach was replaced by Wu's {\it how mathematics} approach to
become the mainstream of studies in the history of mathematics in China in the last quarter of a century, a new
movement is certain to supersede the old one sooner or later. For a historian of mathematics, after the what's and
the how's have been figured out, the problem {\it why} mathematics was done should be addressed.

Following this topic of {\it why mathematics}, research should shift, to some extent, from mathematics in history
to the history of mathematics. Under these circumstances, plenty of new problems are raised for us. Mathematics in
ancient China and other old civilizations, for instance, will be placed in the context of the whole history of
mathematics. The diversity of mathematics in different civilizations will give us a more distinct picture of the
history of mathematics.

{\bf Acknowledgements.} The author is grateful to Prof. Michio Yano and Prof. Wenlin Li for their invaluable
comments on an earlier draft of this article, and to Mr. John Moffett for help with editing. Research for it was
supported by the Japan Society for the Promotion of Science (JSPS, P00019).

\label{lastpage}

\end{document}